# Statistical tools to assess the reliability of self-organizing maps

Eric de Bodt[1], Marie Cottrell[2], Michel Verleysen[3]

[1] Université Lille 2, ESA, Place Deliot, BP 381, F-59020 Lille, France and

Université catholique de Louvain, IAG-FIN, 1 pl. des Doyens,

B-1348 Louvain-la-Neuve, Belgium

[2] Université Paris I, SAMOS-MATISSE, UMR CNRS 8595

90 rue de Tolbiac, F-75634 Paris Cedex 13, France

[3] Université catholique de Louvain, DICE, 3, place du Levant,

B-1348 Louvain-la-Neuve, Belgium

**Acknowledgements**

Michel Verleysen is a Senior Research Fellow of the Belgian National Fund for Scientific Research (FNRS).

# Statistical tools to assess the reliability of self-organizing maps

Eric de Bodt[1], Marie Cottrell[2], Michel Verleysen[3*]

[1] Université Lille 2, ESA, Place Deliot, BP 381, F-59020 Lille, France and

Université catholique de Louvain, IAG-FIN, 1 pl. des Doyens,

B-1348 Louvain-la-Neuve, Belgium

[2] Université Paris I, SAMOS-MATISSE, UMR CNRS 8595

90 rue de Tolbiac, F-75634 Paris Cedex 13, France

[3] Université catholique de Louvain, DICE, 3, place du Levant,

B-1348 Louvain-la-Neuve, Belgium

**Abstract**

Results of neural network learning are always subject to some variability, due to the sensitivity to initial conditions, to convergence to local minima, and, sometimes more dramatically, to sampling variability. This paper presents a set of tools designed to assess the reliability of the results of Self-Organizing Maps (SOM), i.e. to test on a statistical basis the confidence we can have on the result of a specific SOM. The tools concern the quantization error in a SOM, and the neighborhood relations (both at the level of a specific pair of observations and globally on the map). As a by-product, these measures also allow to assess the adequacy of the number of units chosen in a map. The tools may also be used to measure objectively how the SOM are less sensitive to non-linear optimization problems (local minima, convergence, etc.) than other neural network models.

# 1. Introduction

Neural networks are powerful data analysis tools. Part of their interesting properties comes from their inherent non-linearity, in contrast to classical, linear tools. Nevertheless, the non-linear character of the methods has also its drawbacks: most neural network algorithms rely on the non-linear optimization of a criterion, leading to well-known problems or limitations concerning local minima, speed of convergence, etc.

It is commonly argued that vector quantization methods, and in particular self-organizing maps, are less sensitive to these limitations than other classical neural networks, like multi-layer perceptrons and radial-basis function networks. For this reason, self-organizing maps (SOM) [1] are often used in real applications, but rarely studied on the point of view of their reliability: one usually admits that, with some "proper" choice of convergence parameters (adaptation step and neighborhood), the SOM algorithm converges to an "adequate", or "useful", state.

This paper aims at defining objective criteria that may be used to measure the "reliability" of a SOM in a particular situation [2]. Reliability does not mean to measure the quantization error, nor the so-called "topology preservation" of the map; many such criterions have already been published in the literature. Rather, we are looking to measure *how much confidence* we can have in the result of a SOM. This will be achieved through statistical tools based on the bootstrap methodology. The use of statistical tools makes it possible to obtain an objective measure of the confidence we may have in a specific result, and to apply hypotheses tests on these measures. This is different from an approach based on Monte-Carlo simulations of a SOM on a particular database; the last could give an idea about the variance of a numerical result, but could not be used to test if the results of the SOM are statistically significant.

For example, the topology preservation property of a SOM is often used to project neighboring input vectors on the same or on neighboring units (centroids) in the SOM grid. We are looking to evaluate if this projection on neighboring units is obtained by chance, or if it is really the result of topological properties of the database. The reliability of both the quantization and the topology preservation of SOMs are studied. Furthermore, it will be shown that the reliability of the quantization may be used to assess the number of units needed in a SOM.

The study of reliability relies on the extensive use of the bootstrap method. In the following of this paper, we will first address the conventional quantization and organization criteria (section 2), then show how we use the bootstrap methodology in the context of SOMs (section 3). The main contribution of the paper is the definition of reliability criteria (section 4); the criteria are illustrated on artificial and real databases (section 5) before giving conclusions.

## 2. Quantization and organization errors in the SOM

### 2.1. Reliability measures

This section introduces the problem of the measures of quantization and organization in SOMs. First, it should be insisted on the fact that the purpose of this paper is not to define new such criterions. Rather, the question we would like to answer is to know if we can trust the result of the SOM convergence on a particular database (both in terms of quantization and organization).

Trusting the results of a SOM, or measuring its "reliability", does not mean that we expect the same result for different SOMs (even of the same database, but for example sampled differently) run on specific data. Indeed, it is well known that a 1-dimensional SOM (string) on a 2-dimensional distribution, or a 2-dimensional map on a 3-dimensional distribution, will have "folds" because of the string or map trying to estimate the dataset correctly. Even if the dimension of the map and the intrinsic dimensionality of data match, different SOMs will lead to different centroid locations, and even different neighborhood properties, due to twists, butterfly effects, rotation and mirroring of the map, etc. For these reasons, as most of the "neural network" algorithms, SOMs are sometimes criticized because one cannot be sure of the significance of the results.

If we carefully examine what are the causes of these variations in the results, we can see two main sources. First, the initialization of the centroids influences their final positions. Secondly, the final result depends on the sampling of inputs used for learning. .

In this paper, we do not consider specifically the problem of variations of results due to the initialization of the SOM. We have shown in previous works [2], based on Monte-Carlo simulations, that the SOM algorithm is highly insensitive to the choice of initial values. On finite size databases, we have experimentally found that variations due to sampling are much larger then those due to different initial conditions of the algorithm.

We are therefore interested in assessing how much the results (quantization, organization) depend on the specific sample used for training the SOM; in other words, we are interested in the sampling distribution of the performance statistics after learning. As most often in applications, we do not know the exact distribution; moreover, the dataset is of limited (and frequently rather small) size. This suggests using a non parametric bootstrap approach in order to generate artificial ("bootstrapped") samples from the original dataset, and make repeated simulations possible in order to estimate these statistics.

Other resampling or subsampling procedures could be used as well (see for example [4]) and are appropriate modifications of standard bootstrap when this one is shown to be inconsistent. As we do not address in this paper the proof of consistency of the proposed approach, we have no specific reasons to use one of those alternative approaches. The study of consistency will be an interesting development of the current work.

## 2.2. Quantization error

One of the two main goals of Self-Organizing Maps is to quantize the data space into a finite number of so-called *centroids* (or *code vectors*). Vector quantization is used in many areas to compress data over transmission links, to remove noise, etc. The squared distance between an observed data $x^j$ and its corresponding (nearest) centroid is the (quadratic) *quantization error*. Summing this quantization error over all data leads to the *distortion* or *intra-class sum of squares* (which are different names for the same error, used respectively in the information theory domain and by statisticians):

$$SSIntra\left(G^1,\text{K},G^U\right)=\sum_{i=1}^{U}\ \sum_{x^j\in V^i}d^2\left(x^j,G^i\right)=\sum_{x^j}\min_{1\le i\le U}d^2\left(x^j,G^i\right) \qquad (1)$$

where $U$ is the number of units (centroids) in the SOM, $G^i$ is the $i$-th centroid, $d$ is the classical Euclidean distance, $V^i$ is the Voronoi region associated to $G^i$, i.e. the region of the space nearer to $G^i$ than to any other centroid, and the sums on $x^j$ cover all observed data.. Usually, the $d$ distance is the standard Euclidean distance. For the sake of simplicity, we will omit in the following of this paper the explicit mention of the dependency of *SSIntra* on the centroid locations $G^1...G^U$; nevertheless, it should never be forgotten that this dependency exists, and furthermore that the centroid locations are changing during the course of the SOM algorithm. Note that the *SSIntra* criterion is related to the often-used

MSE (Mean Square Error) by a constant factor equal to the number of observations. Expression (1) is the empirical estimation of the theoretical distortion

$$TDIST(\,G^1,\Lambda\,,G^U\,) = \sum_{i=1}^{U} \int_{x\in V_i} d^2\big(x,G^i\big) p\,(dx) = \int_X min_{1\le i\le U}\ d^2\big(x,G^i\big) p\,(dx) \qquad (2)$$

where $p$ is the distribution of the input and $X$ is the input space.

*TDIST* is a very complex function of the distribution $p$ since $G^1$,...., $G^U$ are the centroids computed as limit points of the stochastic algorithm and depend also on the distribution $p$. *SSIntra* is a statistic (calculable from the data) but depending on the unknown distribution $p$, and is precisely what we want to evaluate here.

Note also that the objective function associated with the SOM algorithm for a constant size of neighborhood and finite data set is *the sum of squares intra-classes extended to the neighbor classes,* (see [3]). But actually, one usually ends with no neighbor for the last iterations of the SOM algorithm; at the end of its convergence, the SOM algorithm thus exactly minimizes the *SSIntra* function.

### 2.3. Measure of organization

Twists, butterfly effects, rotation and mirroring of the map may lead to different neighborhood properties of SOMs. However, this does not preclude the use of SOMs in many settings, whether the intrinsic dimensionality of data and the dimension of the map match or not. SOM users know that, because of the twists and folds, one can rarely be confident in the fact that two *specific* neighboring data will be projected on the same or neighboring centroids. However, one can usually be confident in the fact that:

1. the probability that two specific neighboring data will be projected on the same or neighboring centroids is high, and

2. in average over the dataset, pairs of neighboring data will be projected on the same or neighboring centroids.

The first case relates to the observation of specific data over several runs of a SOM or several SOMs, while the last corresponds to spatial averaging of results over a SOM.

There exist several criteria aimed to measure the second case. However, they are not appropriate to measure the first one. Moreover, they give no information about the fact that, over several runs of a SOM (with for example different sequences of presentation of data), the *same* pairs of data will be, in average, projected on the same or neighboring centroids. Having a measure of this last question is precisely what is needed to appreciate if the user of a SOM can be *confident* in the fact that what he/she observes (neighborhood relations on *specific* data) on a single run will be repeated (or not) if other runs were made or other SOMs used (in other words, is the observed results on specific data are reliable or not).

Several measures have been proposed in the literature, in order to measure if the resulting SOM preserves the topology or not. By preserving the topology, we mean that observations close in the input space should be projected to close centroids in the SOM. Among these measures, we can quote the following ones (this list is not exhaustive).

- Demartines [5] plots in a *dy-dx diagram* the distance between the centroids (in the SOM) onto which are projected two specific observations, versus the distance between these observations in the input space; each pair of points in the database is considered. When the topology is locally preserved, the left of the resulting diagram is close from a straight line.

- Villmann [6] who detects the topology preservation of a map through a so-called *topographic function* that measures how neighboring reference vectors in the input space are mapped to neighboring centroids on the map

The important point here is to notice that all these criteria are *global* measures, designed to evaluate if the resulting map is approximately unfolded, or, to the contrary, if some twists (like the well-known butterfly effect) occur. They can also be used to evaluate if, globally on the entire database, local topology is preserved. However none of these criteria is able to test if *the fact that two specific observations are projected on neighboring centroids on the SOM is meaningful or not*. The *STAB* criterion and its histogram as defined further in this paper are an attempt to answer this question.

### 3. Bootstrap

The main idea of the classical bootstrap [7-8-9] is to use the so-called "plug-in principle". Let $F$ be a probability distribution depending on an unknown parameter vector $\theta$. Let $\boldsymbol{x}=x^1, x^2,..., x^n$ be the observed sample of data and $\hat{\boldsymbol{\theta}} = T(\boldsymbol{x})$ an estimate of $\theta$. The bootstrap consists in using artificial samples (called *bootstrapped samples*) with the same empirical distribution as the initial data set in order to guess the distribution of $\hat{\boldsymbol{\theta}}$. Each *bootstrapped sample* consists in $n$ uniform drawings with replacements from the initial sample. If $x*$ is a bootstrapped sample, $T(x*)$ will be a bootstrap replicate of $\hat{\boldsymbol{\theta}}$.

This main idea of the bootstrap may be declined in several ways. In particular, when the evaluation of a statistic $V(x)$ requires non-linear optimization, the well-known problems, or limitations, related to local minima and convergence are met. It may thus happen that different local minima are reached when $V(x*)$ is evaluated for different bootstrapped samples. To avoid this problem, various bootstrap settings can be used. We can speak about [10]:

- Common Bootstrap (**CB**) when each evaluation of $V(x)$ is initialized at random;

- Local Bootstrap (**LB**) when the initial values of each evaluation are kept fixed;

- Local Perturbed Bootstrap (**LPB**) when a small perturbation is applied to the initial conditions obtained as with the Local Bootstrap.

As our purpose is to examine the variability (or the sampling distribution) of some parameters when they are evaluated through different (bootstrapped) samples, but keeping all other conditions (including initial conditions) unchanged, we will use Local Bootstrap. Moreover, as mentioned above (section 2.1), previous works [2] based on Monte-Carlo simulations have shown that the SOM algorithm is highly insensitive to initial conditions.

In our case, the statistics to evaluate are the distribution of SOM quality criteria (both quantization and organization). We know that *SSIntra* is a measure of the quality of quantization; nevertheless, when running a single SOM on a dataset, one cannot be sure if the resulting *SSIntra* value is reliable or not. In our case, it is not possible to have an a priori about the shape of the *SSIntra* distribution. Indeed, having a finite-size database and without knowledge of the exact distribution of observations $x^i$, it is impossible to establish the exact sampling distribution of *SSIntra*. The non parametric bootstrap is therefore used to estimate empirically the distribution of *SSIntra* and, on this basis, to make it possible to evaluate the level of confidence in *SSIntra*.

## 4. Reliability criteria

In this section, we will define original criteria to assess the reliability of the conclusions drawn from the convergence of a SOM on a specific database. Section 4.1 examines the stability of the quantization in a SOM, and how this concept can be taken into consideration to estimate the number of units necessary in the map. Section 4.2 defines the concept of neighborhood relation between two specific observations, and how to use this concept to evaluate the reliability of the topology preservation.

### 4.1. Stability of the quantization in the SOM

The question here is to evaluate the stability of the quantization error in a SOM. In that purpose, we define the relative error, i.e. the *coefficient of variation,* of the *SSIntra* quantization error, as

$$CV(SSIntra) = 100 \frac{\sigma_{SSIntra}}{\mu_{SSIntra}} . \qquad (3)$$

The mean $\mu_{SSIntra}$ and the standard deviation $\sigma_{SSIntra}$ are calculated on the different values of the *SSIntra* error obtained for each bootstrapped case. Obviously, a small *CV* means that the different values of *SSIntra* are close one from another, and therefore that we may have confidence in these values. Note that when repeating the convergence of the SOM, we are not looking to the location of the centroids themselves, but on how they quantify the space in average.

As it will be shown through experimental results in section 5, drawing $CV(SSIntra)$ with respect to the number of units in the SOM can help to assess an adequate number of units in the map. Indeed a sudden increase in the curve usually reveals the fact that some centroids switch from one cluster of the input distribution to another, in other words that there is some instability in the placement of the centroids. A natural consequence is thus to choose a number of centroids immediately below the location of this increase in the curve. This kind of graphical method is very empirical, but very useful and used in other areas, as for instance to choose the number of components to be kept after a Principal Component Analysis, or the number of classes in an ascending hierarchical classification.

The number *B* of bootstrap replications, needed for a sufficiently accurate evaluation of $CV(SSIntra)$, has not been discussed so far. Efron [8] suggests to take between *B*=50 and *B*=200. Our experiments were made with *B*=100; we did not find any improvement in taking a larger value. An example will illustrate this choice in section 5.

**4.2. Stability of the neighborhood relations in the SOM**

**Stability of the neighborhood for a specific pair of observations**

Besides quantization, the second main goal of the SOM is the so-called *topology preservation*, which means that close data in the input space will be quantized by either the same centroid, either two centroids that are close one from another on a predefined string or grid. Often, for example when the SOM is used as a visualization tool, it is desirable to have an objective measure of this *neighborhood* property. We first define

$$NEIGH_{i,j}^{b}(r) = \begin{cases} 0 & if\ x^i\ and\ x^j\ are\ not\ neighbor\ within\ radius\ r \\ 1 & if\ x^i\ and\ x^j\ are\ neighbor\ within\ radius\ r \end{cases} \qquad (4)$$

Being *neighbors within radius r* means that the two observations $x^i$ and $x^j$ are projected on two centroids in the SOM, the distance between these centroids being smaller or equal to $r$. If the radius $r$ is 0, it means that we evaluate if the two data are projected on the same centroid; if $r = 1$, it means that we evaluate if the two data are projected on the same centroid *or* on the immediate neighboring centroids on the string or grid (2 on the string, 8 on the grid), etc. Superscript $b$ in equation (4) means that this result is obtained on the bootstrap sample $b$; of course it may happen that two observations $x^i$ and $x^j$ are projected on neighboring centroids after one simulation and not after another; this is exactly what we are looking to evaluate.

Figure 1 shows a simple example of a 4x5 rectangular Kohonen map, where four observations $x^1$ to $x^4$ are projected; the four observations are mentioned in the box corresponding to their respective nearest centroid. Table 1 shows a few examples of $NEIGH_{i,j}(r)$ calculated from this map; the infinite-norm distance (i.e. the maximum of the distances in each direction) is used to measure the distance $r$ between centroids on the map.

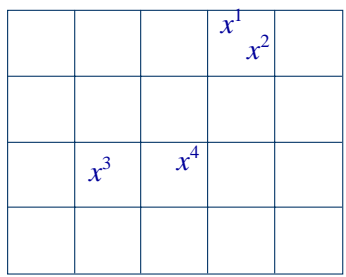

Figure 1. Four observations projected on a 4x5 rectangular SOM.

Table 1. Some examples of *NEIGH* measures on the map illustrated in Figure 1.

| $NEIGH_{1,2}(0) = 1$ | $NEIGH_{3,4}(0) = 0$ | $NEIGH_{1,3}(0) = 0$ |
|---|---|---|
| $NEIGH_{1,2}(1) = 1$ | $NEIGH_{3,4}(1) = 1$ | $NEIGH_{1,3}(1) = 0$ |
| $NEIGH_{1,2}(2) = 1$ | $NEIGH_{3,4}(2) = 1$ | $NEIGH_{1,3}(2) = 1$ |
| $NEIGH_{1,2}(3) = 1$ | $NEIGH_{3,4}(3) = 1$ | $NEIGH_{1,3}(3) = 1$ |
| $NEIGH_{1,2}(4) = 1$ | $NEIGH_{3,4}(4) = 1$ | $NEIGH_{1,3}(4) = 1$ |

We then define $STAB_{i,j}$ as the average of $NEIGH_{i,j}$ over all bootstrap samples:

$$STAB_{i,j}(r) = \frac{\sum_{b=1}^{B} NEIGH_{i,j}^{b}(r)}{B},$$  (5)

where $B$ is the total number of bootstrapped samples.

Still considering the two observations $x^i$ and $x^j$, we are looking to see if *the fact that these two specific observations are projected on neighboring centroids on the SOM is meaningful or not.* One way to evaluate this stability is to see if the $NEIGH_{i,j}(r)$ measure always, or often, takes the same value (0 or 1) over all bootstrap samples. A perfect stability would thus lead $STAB_{i,j}(r)$ to be always 0 ($x^i$ and $x^j$ are never neighbors) or 1 (they are always neighbors).

**Significance of the neighborhood for a specific pair of observations**

The further step is to study the significance of the $STAB_{i,j}(r)$ statistics, by comparing it to the value it would have if the observations fell in the same class (or in two classes distant of less than $r$) in a completely random way (unorganized map). Unorganized maps are taken here as a reference for a *specific* pair of observations $i,j$; indeed if this specific pair of observations only falls in the same class (or in two classes distant of less than $r$) *by chance*, then the $STAB_{i,j}(r)$ statistic will be approximately equal to the value it would take in unorganized maps. Comparing the $STAB_{i,j}(r)$ values in these two situations (organized and unorganized maps) is thus a way to make possible the use of a conventional statistical test to check if the $STAB_{i,j}(r)$ statistic in the organized case is significant or not.

Let $U$ be the total number of classes and $v$ the size of the considered neighborhood. The size $v$ of the neighborhood can be computed from the radius $r$ by $v = (2r + 1)$ for a one-dimensional SOM map (a string); and $v = (2r + 1)^2$ for a two-dimensional SOM map (a grid), if edge effects are not taken into account, or in other words if $v$ is small compared to $U$; this assumption is reasonable, as there is no reason to study the significance of a neighborhood relation if the neighborhood covers a too large part of the map. For a *fixed pair of observations $x^i$ and $x^j$*, with random drawings, the probability of being neighbor in the random case is therefore $v/U$, (it is the probability for $x^j$ to be neighbor of $x^i$ by chance, once the class in which $x^i$ falls is determined).

If we define a Bernoulli random variable with probability of success $v/U$, (where success means: "$x^i$ and $x^j$ are neighbors"), the number $Y$ of successes on $B$ independent bootstrapped trials is distributed as a Binomial distribution, with parameters $B$ and $v/U$. Therefore, it is possible to build a test of the hypothesis $H_0$ "$x^i$ and $x^j$ are only randomly neighbors" against the hypothesis $H_1$ "the fact that $x^i$ and $x^j$ are neighbors or are not neighbors is meaningful".

In order to approximate the binomial random variable by a Gaussian variable, making the hypothesis test easier, we have to check the classical conditions: $B$ has to be large enough (i.e. greater than 30, which is true since we took $B$=100), $B\,v/U > 10$ and $B(1 - v/U) > 10$). For example, for $U$=49, $r$=1, $B$ has to be greater than 50 for the first condition, and greater than 13 for the second one. For $U$=49, $r$=0, the conditions lead to $B > 500$, what becomes too large: the exact binomial distribution of $Y$ must be used. Note that with our $B$=100 choice, the above conditions hold when $v > U/10$ and $v < 9U/10$, corresponding to a neighborhood size that is neither too small not too large; a too large neighborhood is not realistic (it would mean to check the significance of a neighborhood relation when the neighborhood covers most of the map...), while a too small one leads to the necessity of using the original Binomial distribution instead of the approximate Gaussian one.

Next, we can build the rejection region of the $H_0$ against $H_1$ test; for example if the Gaussian approximation is valid, we conclude with a test level of 5% to $H_1$ if $Y$ is less than

$$B\,\frac{v}{U} - 1.96\sqrt{B\,\frac{v}{U}\left(1 - \frac{v}{U}\right)} \text{ or greater than } B\,\frac{v}{U} + 1.96\sqrt{B\,\frac{v}{U}\left(1 - \frac{v}{U}\right)}.$$

It must be mentioned that, because of the random drawing process in the bootstrap, $B$ depends on the pair $(x^i, x^j)$. Indeed the $NEIGH_{i,j}(r)$ only exists if the bootstrap sample contains both observations $x^i$ and $x^j$. We follow the same approach as in [7], which consists in evaluating $STAB_{i,j}(r)$ only on the samples that contain both observations $x^i$ and $x^j$.

**Histogram of the stabilities over all pairs of observations**

Finally, while the $STAB_{i,j}(r)$ is designed to examine a single pair of observations $x^i$ and $x^j$, it is possible to get an idea of the *average stability of the map* by drawing an histogram of the $STAB_{i,j}(r)$ indicator taken over all pairs of observations (the value of $r$ must still be fixed in advance). As mentioned above, a perfect stability would lead $STAB_{i,j}(r)$ to be always 0 or 1; a perfectly stable map would thus lead to an histogram with a first peak at 0, a second peak at 1, and a flat zero curve between these two values. Note that even in such ideal case, the two peaks will not have equal heights; for example when the map is large and the value of $r$ is small, the number of not neighbor pairs is much larger than the number of neighbor ones, leading to a higher peak at 0 in the histogram.

To assess the histogram obtained on a real –not ideally organized– map, one should compare it to the histogram obtained on an unorganized map: in this case there is no reason that neighboring input vectors will be projected on neighboring centroids in such a map, except by chance.

The SOM is organized if its $STAB_{i,j}(r)$ histogram is far from the histogram of an unorganized map. The same arguments as above lead to the conclusion that the distribution of $STAB_{i,j}(r)$ in an unorganized map will follow a Binomial distribution. Nevertheless, the parameters of this Binomial distribution should be chosen carefully for a faithful comparison:

- The probability of success $v/U$ should be modified due to the edge effects on the map: only the centroids in the center of the map have $v = (2r + 1)^2$ neighbors, those on the edges having less neighbors. This modification may be neglected for large maps and small $r$, but not in other cases. As an exemple, $v/U$ should be replaced by 0.15 for a 7x7 map with $r = 1$.

- The number of trials $B$ depends on the pair of observations. The value of $B$ considered for the histogram should thus be the mean of all values of $B$ obtained after simulation.

- As a Binomial distribution takes values between 0 and $B$, the $x$-axis of its histogram must be resized by a factor $B$ in order to be compared to the histogram of the $STAB_{i,j}(r)$ indicator.

In the simulation section, we will use the cumulative histograms for a reliable comparison.

### 4.3. Computational complexity

As with all procedures involving bootstrap and resampling, the computational complexity cannot be ignored. In particular, the number of operations needed to compute the whole set of *STAB* indicators, with $N$ samples in the dataset, $U$ centroids in the map, and $B$ bootstrapped samples, is in O($NUB + N^2B$) (the first term accounts for searching the winning unit for each data in a bootstrapped sample, and the second term for checking if each pair of data is projected on neighboring centroids or not); moreover, it is necessary to train $B$ SOMs. A few comments can however be done:

- if the $STAB_{i,j}(r)$ value is looked for only a *specific* pair of observations $x^i$, $x^j$, the complexity reduces to O($UB$) (training $B$ SOMs remains necessary);

- if the number $N$ of samples in the dataset is large, and if only the histograms of the *STAB* indicators are needed (no information about a *specific* pair of observations), then it is not needed to use all $N$ samples to estimate these histograms; the number of samples used should remain much larger than $U$ though;

- more generally, the *STAB* indicators and histograms may be used to assess the reliability of SOMs in a class of situations (characterized by the complexity and intrinsic dimensionality of data, the size of the dataset, the size of the map, etc.) through a limited number of trials, and extending the conclusions to similar situations; though approximate, this way of working can give some insights about the question whether or not it is necessary, for a specific map on a specific database, to calculate in depth the measures presented in this paper, in order to have an idea about the reliability of the results.

### 5. Experiments

We describe in this section a number of experiments that have been carried out on real and artificial databases, in order to illustrate the concepts introduced in section 4.

All databases used for these simulations are available on the Web site [11]. Simulations have been carried out with 1-dimensional SOM strings or 2-dimensional SOM grids. Results have been obtained with the LB (Local Bootstrap) method.

### 5.1. Real database: macroeconomic situation of countries

The POP_96 database contains seven ratios measured in 1996 on the macroeconomic situation of 96 countries: annual population growth, mortality rate, analphabetism rate, population proportion in high school, GDP per head, GDP growth rate and inflation rate. This dataset has been first used in [12] in the context of data analysis by SOMs.

**Stability of the quantization**

Figure 2 shows the $CV(SSIntra)$ values obtained on the POP_96 database, with increasing number of centroids used in a two-dimensional map.

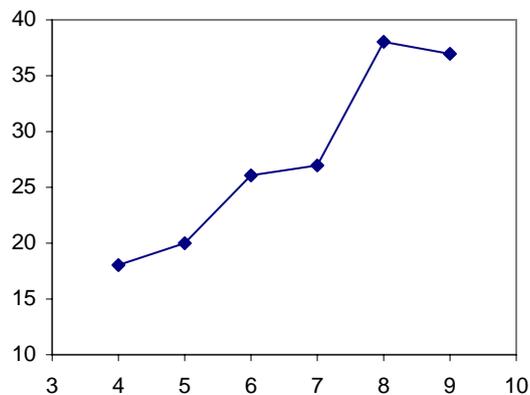

Figure 2. Evolution of $CV(SSIntra)$ with increasing numbers of centroids in a two-dimensional SOM. The *x*-axis shows the number of units in each direction of the map (4x4 to 9x9 units).

As expected, $CV(SSIntra)$ increases with the size of the map: a larger variability in the quantization error is observed for large maps, due to the fact that input vectors switch from one Voronoi region to the other between simulations (a Voronoi region is the portion of the space nearer from one centroid than from any other one).

Nevertheless, it is clear in Figure 2 that a larger increase in $CV(SSIntra)$ is observed at the transition between a 7x7 and a 8x8 map. In order to use a sufficiently large map and, at the same time, to limit the variability of the results, we choose from Figure 2 to use a 7x7 SOM on the POP_96 database.

As mentioned above, a number $B$=100 of bootstrap replications has been chosen in our simulations. Figure 3 illustrates this choice, in the particular case of $CV(SSIntra)$, on the POP_96 database. The value of $CV(SSIntra)$ has been computed 30 times, for each $B$ between 20 and 200 by steps of 20, the purpose being to evaluate if different estimates of $CV(SSIntra)$ (for a specific value of $B$) will lead to similar values. If this is the case, the specific value of $B$ may be considered as large enough; if not, it must be increased. Figure 3 shows the standard deviation of $CV(SSIntra)$ estimated on 30 trials, for each value of $B$; it clearly shows that this standard deviation decreases with $B$ until $B=100$, and remains approximately constant for larger values, justifying the choice. Similar behaviors were observed on other databases and with other reliability criteria than $CV(SSIntra)$.

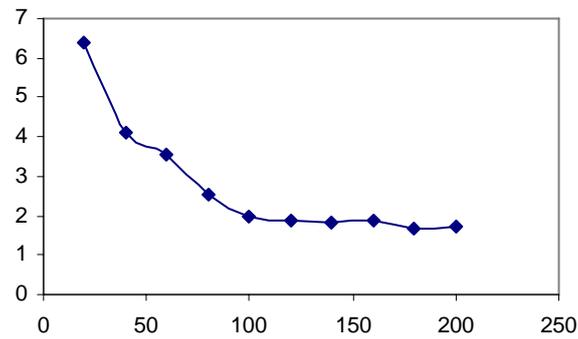

Figure 3. $\sigma(CV(SSIntra))$ as a function of the number $B$ of bootstrap replications.

**Stability of the neighborhood relations**

Table 2 shows the $STAB_{i,j}(r)$ values obtained on a subset of pairs of countries in the POP_96 database. The values marked with * are those which are significant at a test level of 5% ($*^{-0}$ if the result is close from 0 –countries are not neighbors on the map– and $*^{-1}$ if they are close from 1 –they are neighbors–). A 7x7 square map has been used for all simulations.

Table 2. Some examples of *NEIGH* measures on POP database.

| $STAB_{i,j}(r)$ / pairs of countries | $r = 0$ | $r = 1$ | $r = 2$ |
|---|---|---|---|
| Israel – Singapore | 0.023 | 0.7 $*^{-1}$ | 0.95 $*^{-1}$ |
| Albania – Ivory Coast | 0 $*^{-0}$ | 0.17 | 0.59 |
| Guyana – Salvador | 0 $*^{-0}$ | 0.57 $*^{-1}$ | 0.86 $*^{-1}$ |
| Finland – Ireland | 0.51 $*^{-1}$ | 0.92 $*^{-1}$ | 1 $*^{-1}$ |

| Namibia – Uruguay | 0 *-0 | 0 *-0 | 0.15 *-0 |
| Bolivia – Tunisia | 0.02 | 0.33 *-1 | 0.73 *-1 |
| South Korea – Israel | 0.03 | 0.31 | 0.72 *-1 |
| Greece – Italy | 0 *-0 | 0.58 *-1 | 0.87 *-1 |

Pairs of countries like Israel – Singapore and Bolivia – Tunisia are easy to interpret: if the neighborhood concept is measured with $r \geq 1$, both pairs are clearly identified as "neighboring" countries, i.e. in this case as countries with similar macroeconomic standards; taking a smaller radius for the neighborhood concept ($r = 0$) does not lead to any conclusion. The Finland – Ireland pair shows the same behavior, reinforced by the fact that the conclusion can already be taken with $r = 0$. The South Korea – Israel pair are identified as neighboring countries too, but only when taking a large neighborhood concept ($r = 2$).

On the contrary, countries like Namibia – Uruguay really have different macroeconomic standards, as they are significantly identified as "non neighbors" even with $r = 2$.

Conclusions on pairs like Guyana – Salvador and Greece – Italy are different. Countries in these pairs are significantly identified as "non-neighbors" when $r = 0$, and as "neighbors" when $r \geq 1$. This result may seem surprising at first sight. Nevertheless, this can be interpreted as the fact that these countries are always *close* on the map, but –always too– not close enough to be projected on the same centroid.

In average, there are 4.7%, 82.4% and 83.1% of pairs with a "significant" neighborhood relation, i.e. that are marked as *-0 or *-1, respectively for $r = 0$, 1 and 2. If the neighborhood relations were at random, the probability that a pair would be marked with a *-0 or *-1 would be 5%, so the average percent of marked pairs would be 5% too. Thus we see that in this case the results are not significant for r=0, but highly significant for r=1 or r=2.

Figure 4 shows the cumulated histograms of the $STAB_{i,j}(r)$ indicator obtained with a 7x7 SOM on the POP_96 database. The three figures a), b) and c) correspond to $r = 0$ to 2 respectively. We used the cumulated histograms rather than the histograms themselves, for the clarity of the figures. The plain line shows the cumulated histogram of organized SOM, the dashed line the cumulated histogram of an unorganized map; the latter is calculated taking into account the comments at the end of Section 4.2. As detailed in Section 4.2, we conclude that a map is significantly organized if the plain and dashed lines are well separated. As expected on this example, the organization is more significant when the radius $r$ considered for the neighborhood relation is larger. Indeed, when taking $r = 0$ for example, one cannot expect the map to be significantly organized; this would mean that two predefined countries would be either always projected on the *same* centroid (as $r = 0$), or never on the same one. As the number of units in the map (49) is of the same order of magnitudes as the number of countries (96), one easily understands that the discretization of the space is too small to expect such result. On the contrary, when taking a larger radius for the neighborhood relation ($r = 1$ means to consider groups of 9 centroids instead of single ones), the significance of the organization increases (Figure 4 b) and c)).

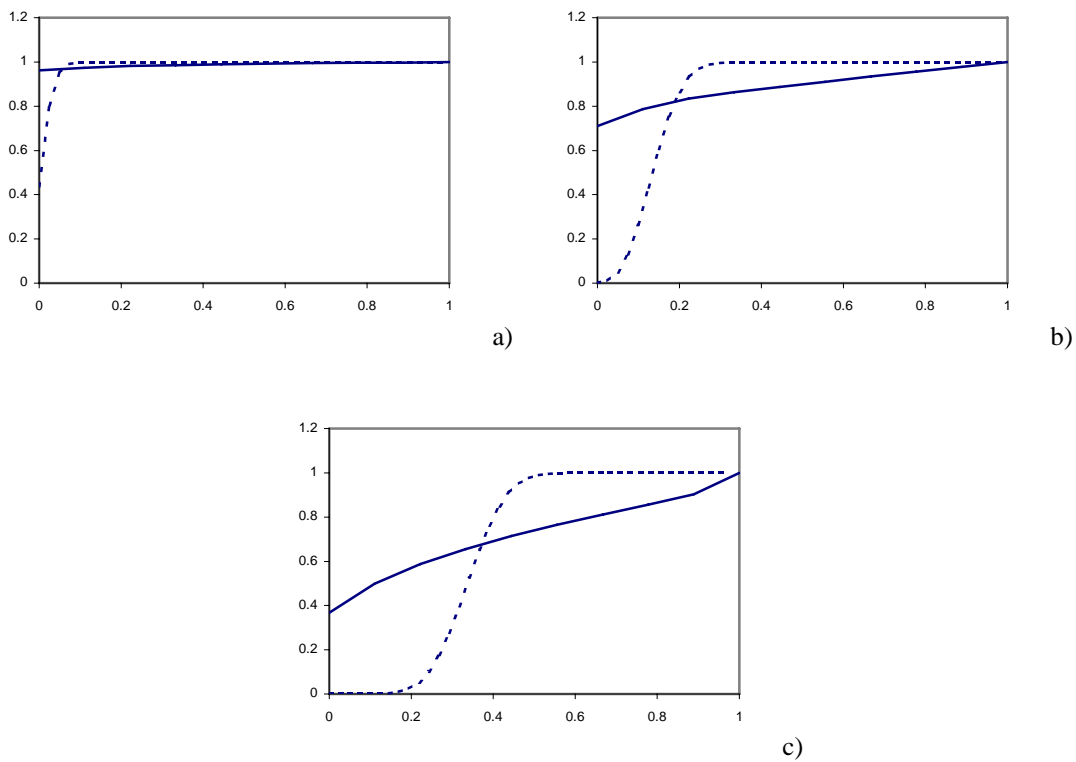

Figure 4. Cumulated histograms of the $STAB_{i;j}(r)$ indicator on the POP_96 database, for a 7x7 SOM. a) $r = 0$; b) $r = 1$; c) $r = 2$. The plain line shows the cumulated histogram of organized SOM, the dashed line the cumulated histogram of an unorganized map. See text for details.

### 5.2. Real database: abalones

The abalone database contains eight numerical indicators measured on more than 4000 abalone shells; indicators include their sizes, weights, age, etc. This database is often used in function approximation benchmarks, where the goal is to predict the age of the abalone as a function of the other variables. In our study we used the eight numerical indicators (including the age) for the SOM learning. A supplementary, non-numerical indicator (the sex of the abalone) contained in the database has not been used. This dataset is available from the UCI repository of machine learning databases [13].

**Stability of the quantization**

A similar study as the one performed on the POP database has been carried out on the abalone one, concerning the evolution of $CV(SSIntra)$ as the size of the SOM increases. Figure 5 shows a rather continuous increase of $CV(SSIntra)$ and therefore, gives no clear indication about the appropriate size of the map. This result is not surprising at the light of the next section. We choose using a 6x6 size in the sequel.

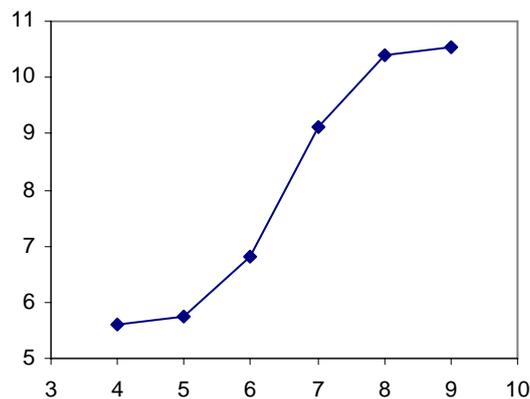

Figure 5. Evolution of $CV(SSIntra)$ with increasing numbers of centroids in a two-dimensional SOM. The *x*-axis shows the number of units in each direction of the map (4x4 to 9x9 units).

**Stability of the neighborhood relations**

Figure 6 shows the cumulated histograms of the $STAB_{ij}(r)$ indicator obtained with a 6x6 SOM on the Abalone database. The three figures a), b) and c) correspond to $r = 0$ to 2 respectively. The plain line shows the cumulated histogram of organized SOM, the dashed line the cumulated histogram of an unorganized map; the latter is calculated taking into account the comments at the end of Section 4.2.

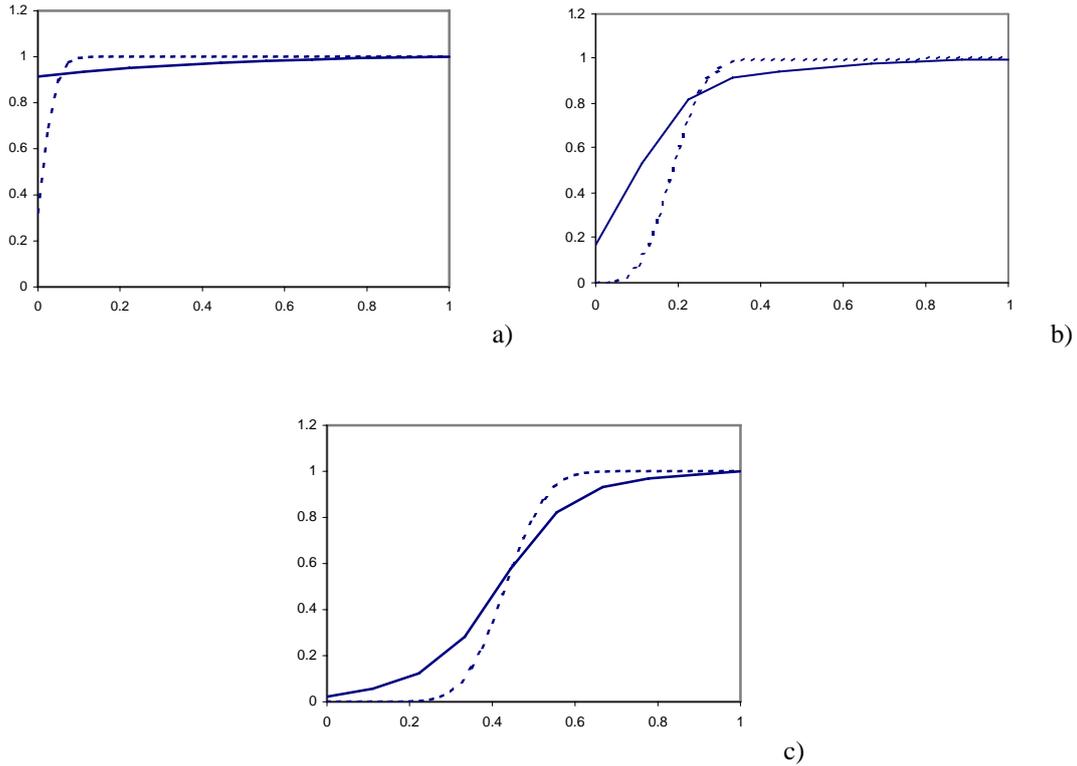

Figure 6. Cumulated histograms of the $STAB_{i;j}(r)$ indicator on the Abalone database, for a 6x6 SOM. a) $r = 0$; b) $r = 1$; c) $r = 2$. The plain line shows the cumulated histogram of organized SOM, the dashed line the cumulated histogram of an unorganized map. See text for details.

We see on the three figures that the organization of the SOM, on the Abalone database, is far from being as good as with the POP database. In all three figures, the plain and dashed lines are close from each other; we may conclude that neighborhood relations between the majority of pairs of data are not significant, i.e. that one cannot have confidence in the fact that, after a specific SOM run, two observations projected on close centroids on the map is not the result of chance. This is due to the nature of the data, which are not well separated into different classes.

### 5.3. Artificial databases: Gaussian distributions

In order to illustrate the concepts on simpler databases, three mono- or multi-modal Gaussian distributions have been used. Gauss_1 contains one cluster of observations. Gauss_2 and Gauss_3 both contain three clusters; the clusters in Gauss_2 have equal variance and some overlap, while those in Gauss_3 have different variance but no overlap. The three distributions are illustrated in Figure 7 a), b) and c) respectively.

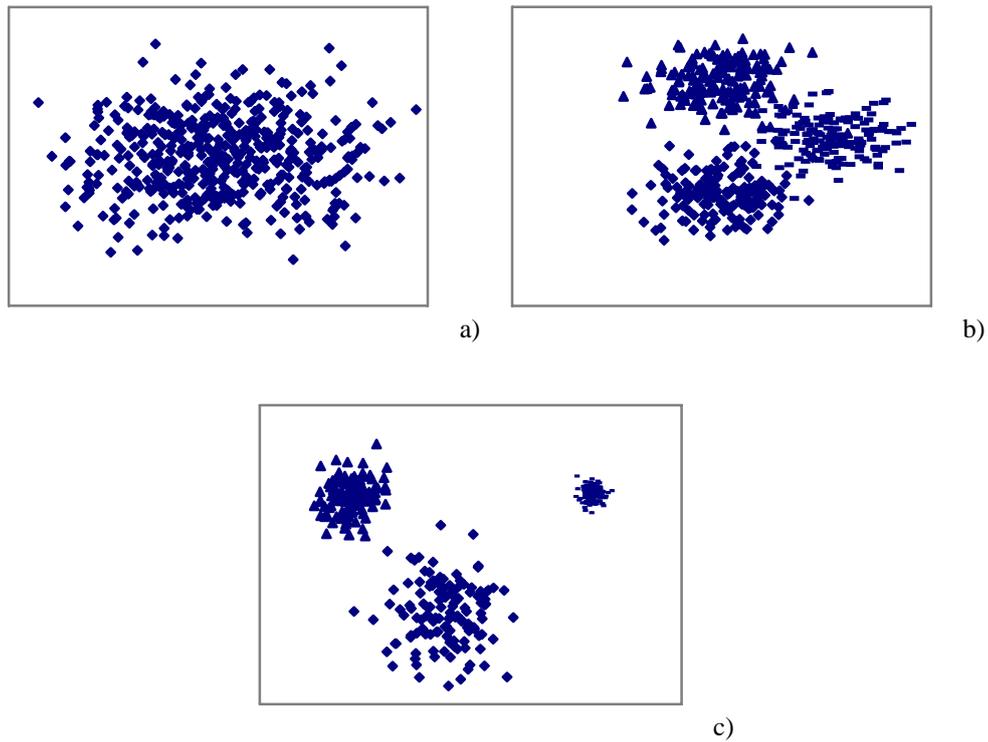

Figure 7.  a) Gauss_1, b) Gauss_2, c) Gauss_3 distributions.

**Stability of the quantization**

Table 3 shows the *CV*(*SSIntra*) values obtained on the three Gaussian databases, with SOM strings of 3, 6 and 9 units respectively.

Table 3.  Coefficient of Variation (*CV*) of the *SSIntra* quantization error, on three Gaussian databases (see text), with one-dimensional SOM strings and variable number of units.

| database / # units | 3 units | 6 units | 9 units |
|---|---|---|---|
| Gauss_1 | 2.04 | 1.85 | 2.40 |
| Gauss_2 | 1.35 | 2.73 | 2.27 |
| Gauss_3 | 1.46 | 8.91 | 5.96 |

An expected result is that the *CV* is low in each case (we remind that the *CV* –Equ. (3)– is expressed as a percentage).  This enforces the idea that the SOM is a reliable method, not falling too much into the traps of local minima encountered with other neural network models (see for example [10] for opposite results obtained with MLP).  We also notice that the *CV* increases with the number of units, in the databases with several clusters (Gauss_2 and Gauss_3). This is due to the fact that some centroids switch from one cluster to another between the different trials, leading to an increase in the variability of the results. We can take advantage of this behavior to have an idea about the number of clusters in a database, and therefore to choose an appropriate number of centroids. Indeed let us have a look to Figure 8, where we detailed the *CV*(*SSIntra*) for the Gauss_1 and Gauss_3 databases.

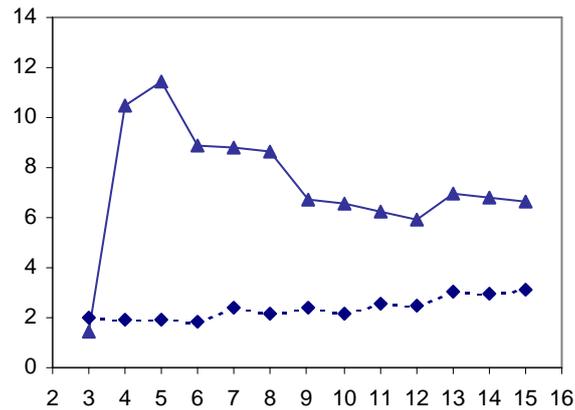

Figure 8. *CV*(*SSIntra*) for the Gauss_1 (squares and dashed line) and Gauss_3 (triangles and plain line) databases, versus the number of centroids in a SOM string.

The plain line in Figure 8 shows the *CV*(*SSIntra*) for the Gauss_3 database, with 3 to 15 units in a SOM string. This line clearly shows that the *CV*(*SSIntra*) strongly decreases at the 5 to 6 units and the 8 to 9 units transition. The reason is that with 5 and 11 units, some centroids switch from one cluster to another between different trials, as already mentioned. However, with 6 and 9 units, the SOM always converges to a situation where there are respectively 2 and 3 units in each of the three clusters. All trials give thus similar results, leading to a lower *CV*(*SSIntra*) value. Diagrams as Figure 8 may be used to evaluate the number of separated clusters in the database (here 3), and to choose an appropriate number of centroids (here a good choice would be a multiple of 3, depending on the expected granularity of the quantization). Of course, this property is not visible in the *CV*(*SSIntra*) for the Gauss_1 database (dashed line in Figure 8), as there is only one cluster in this case.

Since the results concerning the stability of the neighborhood relations are qualitatively similar to those obtained on the previous and the next examples, and since more conclusions can be drawn from the comparison between two databases in the next example, we omit here the results on the neighborhood relations on the Gauss databases.

## 5.4. Artificial databases: 3-dimensional uniform and horseshoe distributions

A third set of databases has been chosen in order to show that two-dimensional SOMs (as commonly used in most situations) perform better (what concerns the quality of the neighborhoods) on databases where the intrinsic dimensionality of data is 2. By intrinsic dimensionality, we mean the effective number of degrees of freedom of the observations.

The database chosen for this example is the well-known horseshoe distribution in a three-dimensional space; its intrinsic dimension is 2, as all observations are situated on a folded 2-surface. The horseshoe distribution is represented in Figure 9. We used a 7x7 SOM grid for the simulations.

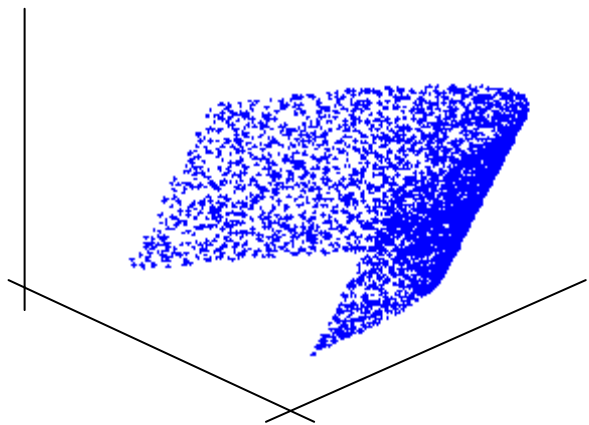

Figure 9. Horseshoe distribution in a 3-dimensional space.

In order to compare the results to a situation where the intrinsic dimensionality of data is larger than 2, we generated a second distribution, where the 3-dimensional points are randomly situated in the whole space. The intrinsic dimensionality of this database is thus 3.

Figure 10 shows the cumulative histograms of *STAB* measures on these two databases (plain lines); dashed lines correspond to the histograms on unorganized maps for comparisons. Figure 10 a) and b) show the results on the uniform distributions, and Figure 10 c) and d) on the horseshoe distribution.

Figures 10 a) and c) have been calculated with a small radius $r = 1$ for the definition of the neighborhood relations, while Figures 10 b) and d) have been calculated with a larger radius $r = 2$. Horizontal axes of all figures have been scaled to 1.

As expected, plain and dashed lines are more different on Figure 10 c) and d) (horseshoe distribution) than they are on Figure 10 a) and b); the difference is more visible between figure 10 b) and d). This means that the neighborhood relations (in particular the larger neighborhoods, i.e. $r = 2$), are better preserved by the SOM on the distribution with intrinsic dimension of 2, than on the distribution with intrinsic dimension of 3. This is nothing else than the phenomenon already mentioned: when the topology of the map does not correspond to the intrinsic dimension of the observations, the SOM has more difficulties to preserve the relations between *all* pairs of points close in the input space. Evaluating the difference between the histograms of *STAB* in organized and unorganized maps is an attempt to measure objectively to what extend the respect of topology is an effective consequence of the convergence of the SOM.

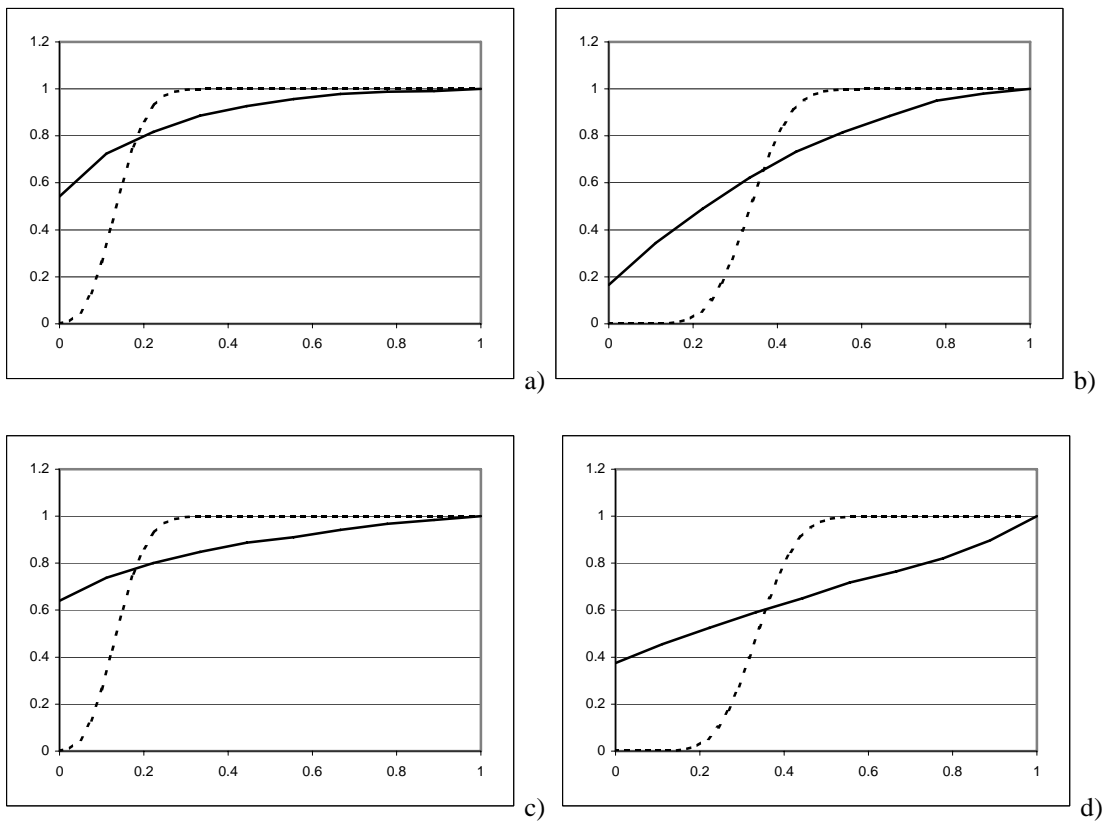

Figure 10. Cumulative histograms of *STAB* measures, on a uniform 3-dimensional database (a) and b)) and on a horseshoe distribution (c) and d)). Plain lines are the result of a 7x7 organized Kohonen map, while dashed lines correspond to unorganized maps. a) and c) have been calculated with $r = 1$, b) and d) with $r = 2$. See text for details.

## 6. Conclusions

This paper describes a set of tools (measures and graphs) to assess the reliability of self-organizing maps. By reliability, we mean the confidence we may have in the result of a specific SOM (after learning) on a specific database. The bootstrap methodology is used as a way to obtain objective tests of statistical significance, that can be used for hypotheses tests in some cases.

Concerning the vector quantization property of SOMs, the reliability is measured by the coefficient of variation of the quantization error. This allows:

- to assess if the value of the quantization error found after the SOM learning is reliable or not, i.e. is reproduced if the same experiment is made again but with different samples;

- as a by-product, to assess if the number of units (centroids) chosen in the map is adequate with respect to the clusters in the data (a bad choice will lead to an sudden increase in the coefficient of variation).

An original result of this paper is to measure the topological property of SOMs separately for each pair of observations, in order to check if the fact that they are neighbors (or not) in a SOM (after learning) is meaningful or not. The bootstrap methodology is used to average over a set of simulations, a measure of the neighborhood relation on a specific pair of observations. This allows:

- to check if the fact that two observations are projected on the same or on neighboring centroids on the SOM is a meaningful result (against the fact that it is a random result);

- by looking to the results of a larger set of pairs, to assess a reasonable size to define the neighborhood concept between two observations.

Concatenating the information about the neighborhood relations on all pairs of observations in a database, in the form of a histogram, and comparing this histogram to the one that would be obtained if the map was not organized, allows:

- to have an overall measure of the organization of the map *and* a point of comparison to assess if this measure is conclusive or not;

- to assess if the dimension of the map (usually dimension 2 – grid – or dimension 1 – string –) is adequate for the database (according to its intrinsic dimension) or if it will lead to a folded map in the input space (for which topology is not optimally preserved).

**Acknowledgements**